\documentclass[final]{siamltex}
\usepackage{amssymb,amsmath,url}
\usepackage[dvipdfm]{graphicx}

\def\pmatrix{ \left( \begin{array} }
\def\endpmatrix{ \end{array} \right) }

\newtheorem{rem}{Remark}
\newcommand{\bx}{\mbox{\boldmath$b$}}
\newcommand{\cx}{\mbox{\boldmath$c$}}
\newcommand{\gx}{\mbox{\boldmath$g$}}
\newcommand{\vx}{\mbox{\boldmath$v$}}
\newcommand{\wx}{\mbox{\boldmath$w$}}
\newcommand{\yx}{\mbox{\boldmath$y$}}
\newcommand{\zx}{\mbox{\boldmath$z$}}
\newcommand{\ox}{\mbox{\boldmath$o$}}
\newcommand{\bfeta}{\mbox{\boldmath$\eta$}}

\newcommand{\hbfy}{\widehat{\yx}}
\newcommand{\hbfz}{\widehat{\zx}}
\newcommand{\hbfw}{\widehat{\wx}}
\newcommand{\cL}{{\mathcal L}}
\newcommand{\cF}{{\mathcal F}}
\newcommand{\cG}{{\mathcal G}}
\newcommand{\cI}{{\mathcal I}}

\def\RR{\mathbb R}
\def\ddt{\frac{\partial}{\partial t}}

\def\urltilda{\kern -.15em\lower .7ex\hbox{\Large \~{}}\kern .04em}
\def\urltilda{{\footnotesize $\sim$}\,}
\title{Parallel Factorizations in Numerical Analysis\thanks{Work developed
within the project ``Numerical methods and software for
differential equations''.}}

\author{Pierluigi Amodio\thanks{Dipartimento di Matematica,
Universit\`a di Bari, Bari, Italy ({\tt amodio@dm.uniba.it}).}
\and Luigi Brugnano\thanks{Dipartimento di Matematica,
Universit\`a di Firenze, Firenze, Italy}  ({\tt
luigi.brugnano@unifi.it}). }

\begin{document}

\maketitle

\begin{abstract}
In this paper we review the parallel solution of sparse linear
systems, usually deriving by the discretization of ODE-IVPs or
ODE-BVPs. The approach is based on the concept of {\em parallel
factorization} of a (block) tridiagonal matrix. This allows to
obtain efficient parallel extensions of many known matrix
factorizations, and to derive, as a by-product, a unifying
approach to the parallel solution of ODEs.
\end{abstract}



\pagestyle{myheadings} \thispagestyle{plain} \markboth{P.\,AMODIO
AND L.\,BRUGNANO}{PARALLEL FACTORIZATIONS IN NUMERICAL ANALYSIS}

\section{Introduction}\label{sec:intro}

The numerical solution of ODEs requires the solution of sparse and
structured linear systems. The parallel solution of these problems
may be obtained in two ways: for BVPs, since the size of the
associated linear system is large, we need to develop parallel
solvers for the obtained linear systems; for IVPs we need to define
appropriate numerical methods that allow to obtain parallelizable
linear systems.

In both cases, the main problem can then be taken back to the
solution of special sparse linear systems, whose solution is here
approached through the use of {\em parallel factorizations},
originally introduced for deriving efficient parallel tridiagonal
solvers \cite{AmBr1,AmBrPo}, and subsequently generalized to block
tridiagonal, Almost Block Diagonal (ABD), and Bordered Almost
Block Diagonal (BABD) systems
\cite{Ametal,AmGlRo1,AmPa1,AmPa2,AmPa3,AmRo}.

With this premise, the structure of the paper is the following: in
Section~\ref{PF} the main facts about {\em parallel
factorizations} and their extensions are briefly recalled; then,
in Section~\ref{par} their application for solving ODE problems is
sketched; finally, in Section~\ref{parareal} we show that this
approach also encompasses the so called ``Parareal'' algorithm,
recently introduced in \cite{LiMaTu01,MaTu02}.

\section{Parallel factorizations}\label{PF}
In this section we consider several parallel algorithms in the
class of partition methods for the solution of linear systems,
\begin{equation} \label{syst}
Ax = f,
\end{equation}
where $A$ is a $n \times n$ sparse and structured matrix, and $x$
and $f$ are vectors of length $n$. We will investigate the
parallel solution of \eqref{syst} on $p$ processors, supposing $p
\ll n$ in order for the number of sequential operations to be much
smaller than that of parallel ones.

\begin{figure}[t]
\begin{center}
\includegraphics*[width=9cm]{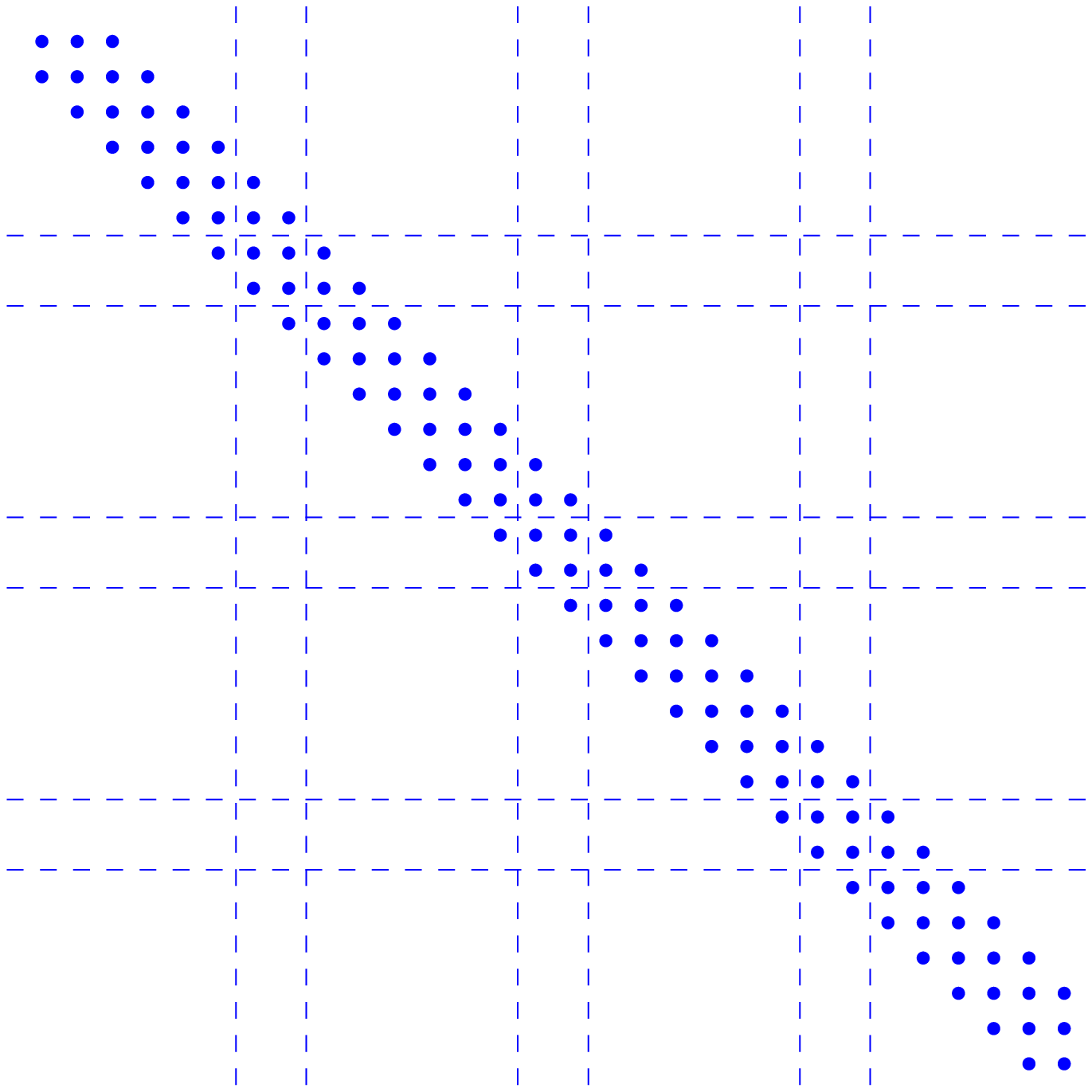}
\vspace*{-1cm}
\end{center}
\caption{\it Partitioning of a banded matrix. Each point
represents a (block) entry of the matrix.} \label{Fig1}
\end{figure}

\begin{figure}[t]
\begin{center}
\includegraphics*[width=9cm]{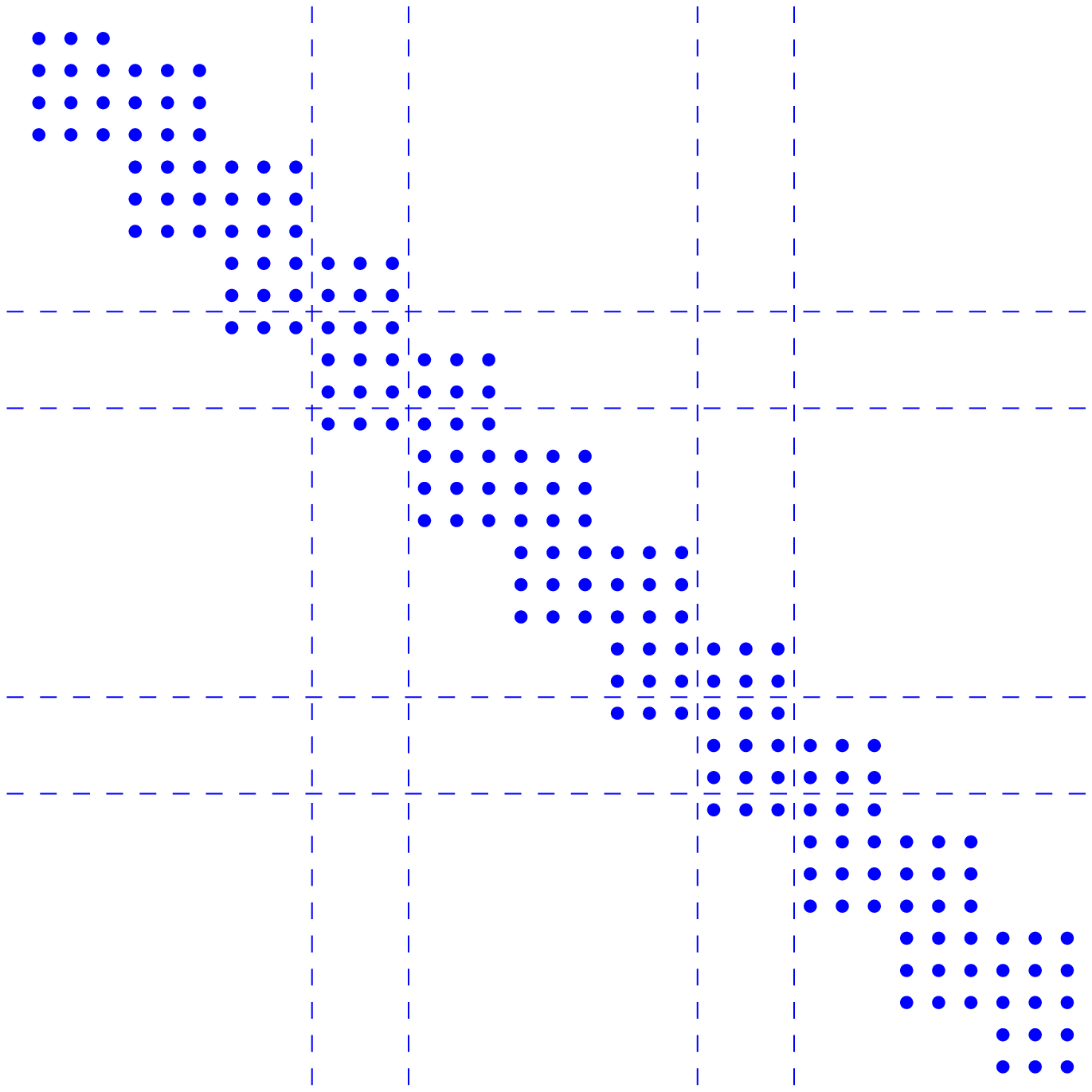}
\vspace*{-1cm}
\end{center}
\caption{\it Partitioning of an ABD matrix. Each point represents an
entry of the matrix.} \label{Fig2}
\end{figure}

The coefficient matrices $A$ here considered are (block) banded,
tridiagonal, bidiagonal, or even Almost Block Diagonal (ABD). All
these structures may be rearranged in the form
\begin{equation} \label{parpart}
A = \pmatrix{ccccccc}
A^{(1)} & \cx_1^{(1)} \\
\bx_1^{(1)^T} & a^{(1)} & \cx_0^{(2)^T} \\
& \bx_0^{(2)} & A^{(2)} & \cx_1^{(2)} \\
&& \bx_1^{(2)^T} & a^{(2)} \\
&&&& \ddots \\
&&&&& a^{(p-1)} & \cx_0^{(p)^T} \\
&&&&& \bx_0^{(p)} & A^{(p)}
\endpmatrix
\end{equation}
where the diagonal blocks are square and the superscript $(i)$
indicates that this block is handled only by processor $i$. The
size of the blocks $a^{(i)}$, $A^{(i)}$, $\bx_j^{(i)}$, and
$\cx_j^{(i)}$ is in general independent of both $i$ and $j$, and
only depends on the sparsity structure of the coefficient matrix
$A$. In particular, the size of the blocks $a^{(i)}$ is quite
important, since the sequential section of the algorithm is
proportional to it. Therefore, the blocks $a^{(i)}$ should be as
small as possible. As an example, if $A$ is (block) tridiagonal,
$a^{(i)}$ reduces to a single (block) entry. Vice versa, in case
of banded (block) matrices, the (block) size of $a^{(i)}$ equals
to $\max(s,r)$, where $s$ and $r$ denote the number of lower and
upper off (block) diagonals (see Figure~\ref{Fig1}), respectively.
In case of ABD matrices, $a^{(i)}$ is a block of size equal to the
number of rows in each block row of the coefficient matrix (see
Figure~\ref{Fig2}). Since row and column permutations inside each
block do not destroy the sparsity structure of the coefficient
matrix, in ABD matrices we may permute the elements inside
$a^{(i)}$ to improve stability properties. Blocks $A^{(i)}$ have
the same sparsity structure as the original matrix, and are
locally handled by using any suitable sequential algorithm.

In order to keep track of any parallel algorithm, we consider the
following factorization \cite{AmBr1,AmBrPo}
\begin{equation} \label{parfact}
A = F\,T\,G,
\end{equation}
where
\begin{equation} \label{parF}
F = \pmatrix{ccccccccc}
N^{(1)} & \ox \\
\vx^{(1)^T} & I & \wx^{(2)^T} \\
& \ox & N^{(2)} & \ox \\
&& \vx^{(2)^T} & I & \wx^{(3)^T} \\
&&& \ox & N^{(3)} & \ox \\
&&&& \vx^{(3)^T} & I \\
&&&&&& \ddots \\
&&&&&&& I & \wx^{(p)^T} \\
&&&&&&& \ox & N^{(p)}
\endpmatrix,
\end{equation}
\begin{equation} \label{parT}
T = \pmatrix{ccccccccc}
\hat I & \ox \\
\ox^T & \alpha^{(1)} & \ox^T & \gamma^{(2)} \\
& \ox & \hat I & \ox \\
& \beta^{(2)} & \ox^T & \alpha^{(2)} & \ox^T & \gamma^{(3)} \\
&&& \ox & \hat I & \ox \\
&&&  \beta^{(3)} & \ox^T & \alpha^{(3)} \\
&&&&&& \ddots \\
&&&&&&& \alpha^{(p-1)} & \ox^T \\
&&&&&&& \ox & \hat I
\endpmatrix,
\end{equation}
\begin{equation} \label{parG}
G = \pmatrix{ccccccccc}
S^{(1)} & \yx^{(1)} \\
\ox^T & I & \ox^T \\
& \zx^{(2)} & S^{(2)} & \yx^{(2)} \\
&& \ox^T & I & \ox^T \\
&&& \zx^{(3)} & S^{(3)} & \yx^{(3)} \\
&&&& \ox^T & I \\
&&&&&& \ddots \\
&&&&&&& I & \ox^T \\
&&&&&&& \zx^{(p)} & S^{(p)}
\endpmatrix,
\end{equation}
$I$, $\hat I$ and $\ox$ are identity and null matrices of
appropriate sizes, and $N^{(i)}S^{(i)}$ is any suitable
factorization of the block $A^{(i)}$. The remaining entries of
$F$, $T$, and $G$ can be derived from \eqref{parfact} by direct
identification.

Factorization  \eqref{parfact} may be computed in parallel on the
$p$ processors. For simplicity, we analyze the factorization of
the sub-matrix identified by the superscript $(i)$ (with obvious
differences for $i=1$ and $i=p$)
\begin{equation}\label{Mi0}
M^{(i)} = \pmatrix{ccc}
0 & \cx_0^{(i)^T} \\
\bx_0^{(i)} & A^{(i)} & \cx_1^{(i)} \\
& \bx_1^{(i)^T} & a^{(i)} \\
\endpmatrix.
\end{equation}
Following \eqref{parpart}--\eqref{parG} we have
\begin{equation} \label{parblock}
M^{(i)} = \pmatrix{ccc}
I & \wx^{(i)^T} \\
\ox & N^{(i)} & \ox \\
& \vx^{(i)^T} & I \\
\endpmatrix
\pmatrix{ccc}
\alpha_1^{(i)} & \ox^T & \gamma^{(i)} \\
\ox & I & \ox \\
\beta^{(i)} & \ox^T & \alpha_2^{(i)} \\
\endpmatrix
\pmatrix{ccc}
I & \ox^T \\
\zx^{(i)} & S^{(i)} & \yx^{(i)} \\
& \ox^T & I \\
\endpmatrix,
\end{equation}
where $$\alpha_2^{(i)}+\alpha_1^{(i+1)} = \alpha^{(i)},$$ and the
other entries are the same as defined in
\eqref{parF}--\eqref{parG}.

Consequently, by considering any given factorization for
$A^{(i)}$, it is possible to derive corresponding parallel
extensions of such factorizations, which cover most of the
parallel algorithms in the class of partition methods. In
particular, the following ones easily derive for matrices with
well conditioned sub-blocks $A^{(i)}$ (this means, for example,
that pivoting is unnecessary or does not destroy the sparsity
structure):

\begin{itemize}

\medskip
\item \textit{$LU$ factorization}, by setting in  \eqref{parblock} $N^{(i)}=L^{(i)}$
and $S^{(i)}=U^{(i)}$, where $L^{(i)}U^{(i)}$ is the $LU$
factorization of the matrix $A^{(i)}$. In this case, the (block)
vectors $\yx^{(i)}$ and $\vx^{(i)}$ maintain the same sparsity
structure as that of $\cx_1^{(i)}$ and $\bx_1^{(i)}$,
respectively, while the vectors $\zx^{(i)}$ and $\wx^{(i)}$ are
non-null fill-in (block) vectors, obtained by solving two
triangular systems.

\medskip
\item \textit{$LUD$ factorization} (which derives from the Gauss-Jordan elimination algorithm),
by setting  in  \eqref{parblock} $S^{(i)}=D^{(i)}$, a diagonal
matrix, and
$$
\pmatrix{ccc}
I & \wx^{(i)^T} \\
\ox & N^{(i)} & \ox \\
& \vx^{(i)^T} & I \\
\endpmatrix =
\pmatrix{ccc}
I & \ox^T \\
\ox & L^{(i)} & \ox \\
& \vx^{(i)^T} & I \\
\endpmatrix
\pmatrix{ccc}
I & \wx^{(i)^T} \\
\ox & U^{(i)} & \ox \\
& \ox^T & I \\
\endpmatrix,
$$
where $L^{(i)}$ and $U^{(i)}$ are lower and upper triangular
matrices, respectively, with unit diagonal. Therefore, $\vx^{(i)}$
and $\wx^{(i)}$ maintain the same sparsity structure as that of
$\bx_1^{(i)}$ and $\cx_0^{(i)}$, respectively, while $\zx^{(i)}$
and $\yx^{(i)}$ are non-null fill-in (block) vectors.

\medskip
\item \textit{cyclic reduction} algorithm \cite{AmBr1,AmBrPo} (see also
\cite{Am,AmMa,AmPa1,AmPa2}), which is one of the best known
parallel algorithms but that, in its original form, requires a
synchronization at each step of reduction. In fact, the idea of
this algorithm is to perform several reductions that, at each
step, halve the size of the system. On the other hand, to obtain a
factorization in the form \eqref{parblock}, it is possible to
consider cyclic reduction as a sequential algorithm to be applied
locally,
$$
M^{(i)} = (\hat P_1^{(i)} \hat L_1^{(i)} \hat P_2^{(i)} \hat
L_2^{(i)} \cdots ) \hat D^{(i)} (\cdots \hat U_2^{(i)} \hat
P_2^{(i)^T} \hat U_1^{(i)} \hat P_1^{(i)^T}),
$$
where $\hat P_i^{(i)}$ are suitable permutation matrices that
maintain the first and last row in the reduced matrix. The
computational cost, which is much higher if the algorithm is
applied to $A$ on a sequential computer, becomes comparable to the
previous local factorizations since this algorithm does not
compute fill-in block vectors. As a consequence, the corresponding
parallel factorization algorithm turns out to be one of the most
effective.

\medskip
\item \textit{Alternate row and column elimination} \cite{Va} which is an
algorithm suitable for ABD matrices. In fact, for such matrices
alternate row and column permutations always guarantee stability
without fill-in. This feature extends to the parallel algorithm,
by taking into account that row permutations between the first
block  row of $A^{(i)}$ and the block containing $\cx_0^{(i)}$
(see (\ref{Mi0})), still make the parallel algorithm stable
without introducing fill-in. Such parallel factorization is
defined by setting $N^{(i)}=P^{(i)}L^{(i)}$ and
$S^{(i)}=U^{(i)}Q^{(i)}$, where $P^{(i)}$ and $Q^{(i)}$ are
permutation matrices and $L^{(i)}$ and $U^{(i)}$, after a suitable
reordering of the rows and of the columns, are $2 \times 2$ block
triangular matrices (see \cite{Ametal} for full details). Finally,
the (block) vectors $\yx^{(i)}$ and $\zx^{(i)}$ maintain the same
sparsity structure as that of $\bx_1^{(i)}$ and $\cx_0^{(i)}$,
respectively, whereas $\wx^{(i)}$ and $\vx^{(i)}$ are fill-in
(block) vectors.

\end{itemize}
\medskip

For what concerns the solution of the systems associated to the
previous parallel factorizations, there is much parallelism
inside. The solution of the systems with the matrices $F$ and $G$
may proceed in parallel on the different processors. Conversely,
the solution of the system with the matrix $T$ requires a
sequential part, consisting in the solution of a {\em reduced
system} with the (block) tridiagonal {\em reduced matrix}
\begin{equation} \label{redmatr}
T_p = \pmatrix{cccc}
\alpha^{(1)} & \gamma^{(2)} \\
\beta^{(2)} & \alpha^{(2)} & \ddots \\
& \ddots & \ddots & \gamma^{(p-1)} \\
&& \beta^{(p-1)} & \alpha^{(p-1)}
\endpmatrix.
\end{equation}
We observe that the (block) size of $T_p$ only depends on $p$ and
is independent of $n$.

\medskip
For a matrix $A$ with singular or ill-conditioned sub-blocks
$A^{(i)}$, the local factorizations may be unstable or even
undefined. Consequently, it is necessary to slightly modify the
factorization \eqref{parblock}, in order to obtain stable parallel
algorithms. The basic idea is that factorization \eqref{parblock}
may produce more than two entries in the {\em reduced system}. In
other words, the factorization of $A^{(i)}$ is stopped when the
considered sub-block is ill-conditioned (or the local
factorization with a singular factor). As a consequence, the size
of the {\em reduced system} is increased as sketched below. Let
then
$$ M^{(i)} = \hat L_1^{(i)} \hat D_1^{(i)} \hat U_1^{(i)},$$
where
$$
\hat L_1^{(i)} = \pmatrix{ccccc}
I & \wx_1^{(i)^T} \\
\ox & N_1^{(i)} & \ox \\
& \vx_1^{(i)^T} & I & \ox^T \\
&& \ox & \hat I & \ox \\
&&& \ox^T & I
\endpmatrix, \quad \hat U_1^{(i)} =
\pmatrix{ccccc}
I & \ox^T \\
\zx_1^{(i)} & S_1^{(i)} & \yx_1^{(i)} \\
& \ox^T & I & \ox^T \\
&& \ox & \hat I & \ox \\
&&& \ox^T & I
\endpmatrix,
$$
$$
\hat D_1^{(i)} = \pmatrix{ccccc}
\alpha_1^{(i)} & \ox^T & \gamma_1^{(i)} \\
\ox & I & \ox \\
\beta_1^{(i)} & \ox^T & \alpha_2^{(i)} & \cx_2^{(i)^T} \\
&& \bx_2^{(i)} & A_2^{(i)} & \cx_3^{(i)} \\
&&& \bx_3^{(i)^T} & \alpha_3^{(i)}
\endpmatrix,
$$
when the sub-block $A_1^{(i)}$ of $A^{(i)}$,
$$
A_1^{(i)} = \pmatrix{cc}
N_1^{(i)}  \\
\vx_1^{(i)^T} & \alpha_2^{(i)}
\endpmatrix
\pmatrix{cc}
S_1^{(i)} & \yx_1^{(i)} \\
& I
\endpmatrix,
$$
is singular, because the block $\alpha_2^{(i)}$ is singular (i.e.,
$\alpha_2^{(i)}=0$, in the scalar case). Then, $\alpha_2^{(i)}$ is
introduced in the {\em reduced system}. By iterating this
procedure on $\hat D_1^{(i)}$, we obtain again the factorization
\eqref{parfact}, with the only difference that now the {\em
reduced matrix} in $T_p$ may be of (block) size larger than $p-1$
(compare with (\ref{redmatr})). However, it may be shown that it
still depends only on $p$, whereas it is independent of $n$
\cite{AmBr2,AmBr3}. Consequently, the scalar section of the whole
algorithm is still negligible, when $n \gg p$.

The parallel algorithms that fall in this class are
\cite{AmBr2,AmBr3}:
\begin{itemize}
\medskip
\item \textit{$LU$ factorization with partial pivoting}, defined by setting $N_1^{(i)} = (P_1^{(i)})^TL_1^{(i)}$ and
$S_1^{(i)}=U_1^{(i)}$ where $P_1^{(i)}$ is a permutation matrix
such that $L_1^{(i)}U_1^{(i)}$ is the $LU$ factorization of
$P_1^{(i)}A_1^{(i)}$. The remaining (block) vectors are defined
similarly as in the case of the $LU$ factorization previously
described.

\medskip
\item \textit{$QR$ factorization}, defined by setting $N_1^{(i)}=Q_1^{(i)}$ and $S_1^{(i)}=R_1^{(i)}$. In this
case both $\wx_1^{(i)}$ and $\zx_1^{(i)}$ are fill-in (block)
vectors while $\vx_1^{(i)}$ and $\yx_1^{(i)}$ maintain the same
sparsity structure as that of the corresponding (block) vectors in
$M^{(i)}$.

\end{itemize}
\medskip

\begin{figure}[t]
\begin{center}
\includegraphics*[width=9cm]{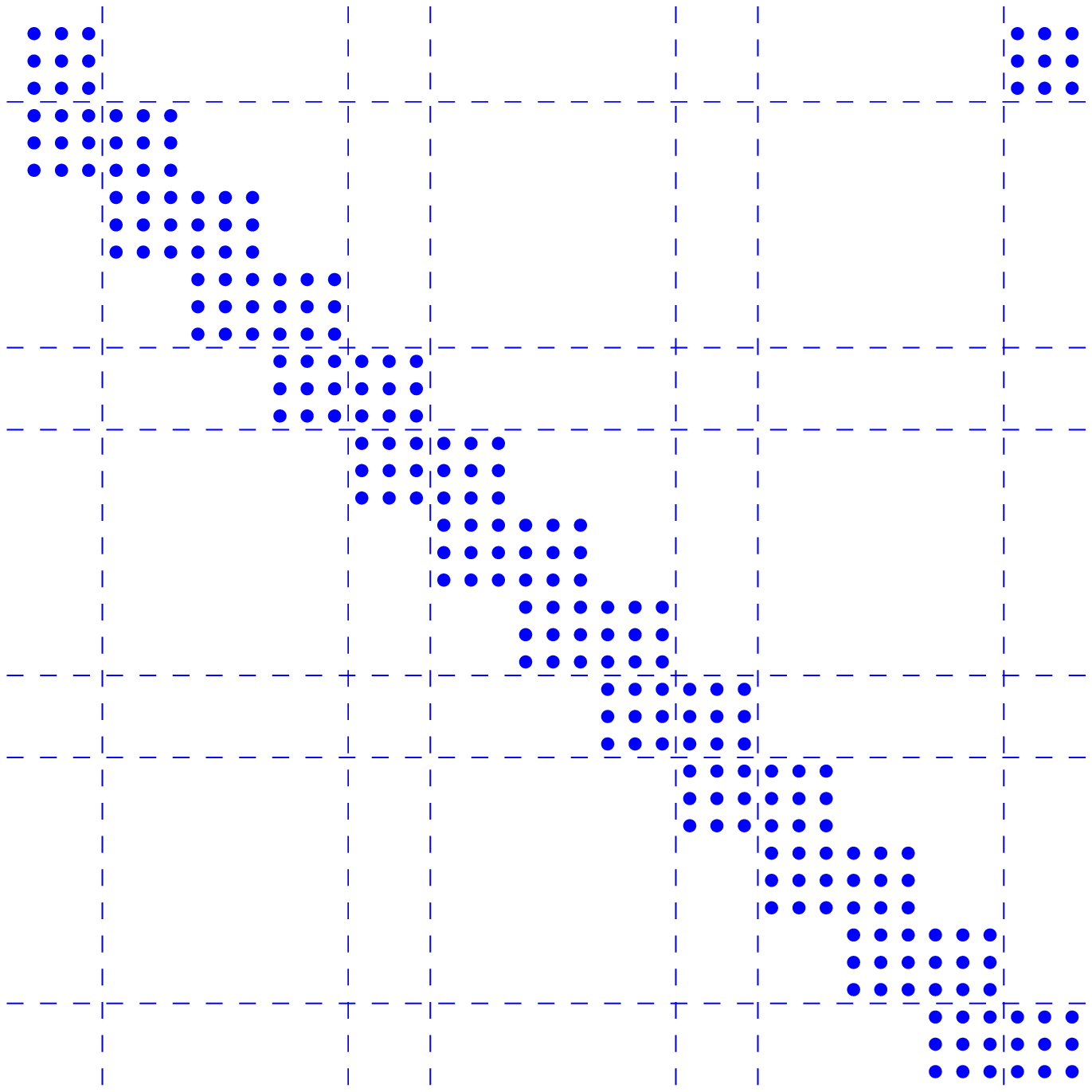}
\vspace*{-1cm}
\end{center}
\caption{\it Partitioning of a BABD matrix. Each point represents
an entry of the matrix.} \label{Fig4}
\end{figure}

\begin{figure}[t]
\begin{center}
\includegraphics*[width=9cm]{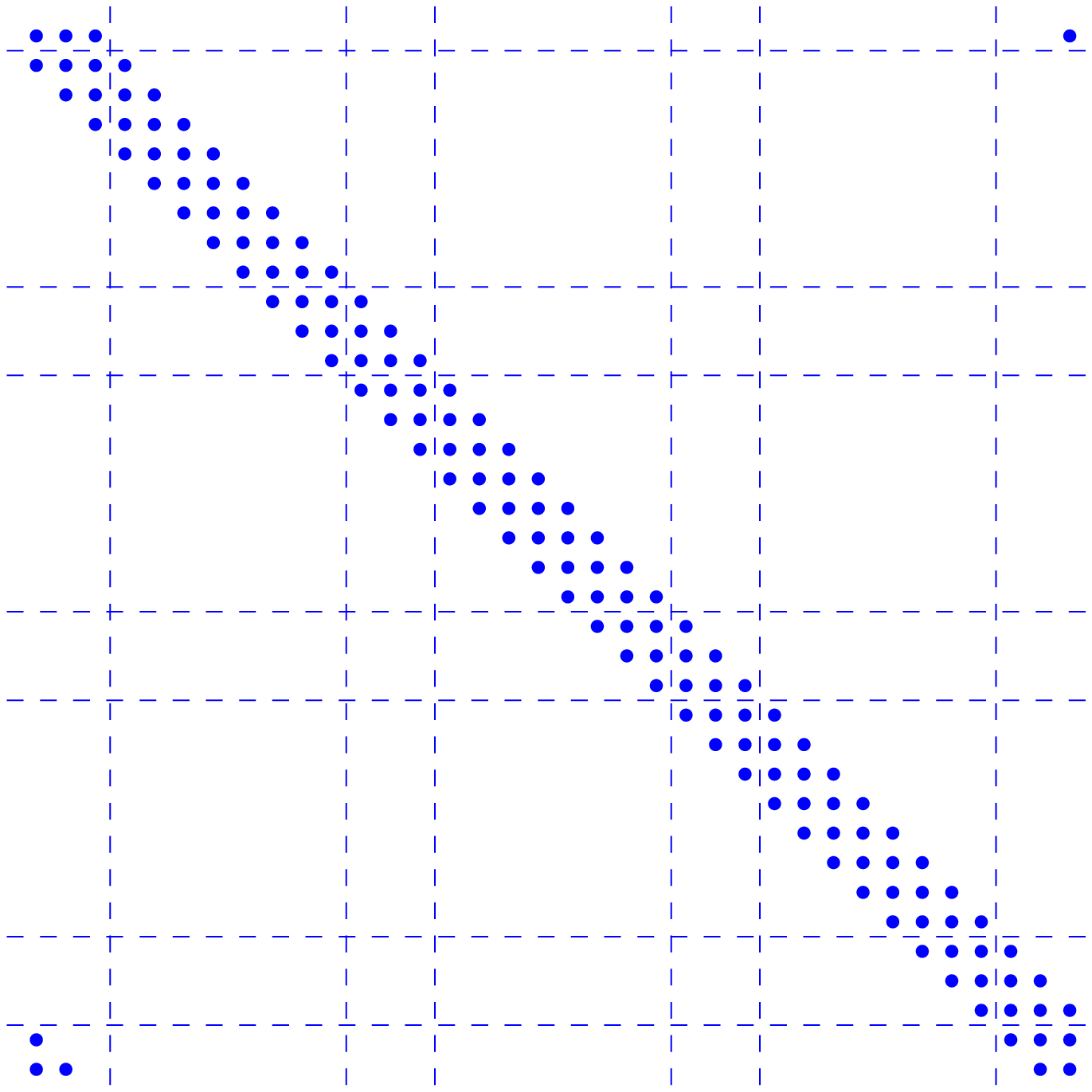}
\vspace*{-1cm}
\end{center}
\caption{\it Partitioning of a matrix with a circulant-like
structure. Each point represents a (block) entry of the matrix.}
\label{Fig3b}
\end{figure}

Factorization (\ref{parfact})--(\ref{parG}), and the corresponding
parallel algorithms mentioned above, are easily generalized to
matrices with additional non-null elements in the right-lower
and/or left-upper corners. This is the case, for example, of
Bordered ABD (BABD) matrices (see Figure \ref{Fig4}) and matrices
with a circulant-like structure (see Figure \ref{Fig3b}).
Supposing the non-null elements are located in the right-upper
corner (this is always possible by means of suitable permutation),
then the coefficient matrix is partitioned in the form
\begin{equation} \label{parpart1}
A = \pmatrix{ccccccccc}
a^{(0)} & \cx_0^{(1)^T} &&&&&&& b \\
\bx_0^{(1)} & A^{(1)} & \cx_1^{(1)} \\
&\bx_1^{(1)^T} & a^{(1)} & \cx_0^{(2)^T} \\
&& \bx_0^{(2)} & A^{(2)} & \cx_1^{(2)} \\
&&& \bx_1^{(2)^T} & a^{(2)} \\
&&&&& \ddots \\
&&&&&& a^{(p-1)} & \cx_0^{(p)^T} \\
&&&&&& \bx_0^{(p)} & A^{(p)} & \cx_1^{(p)} \\
&&&&&&& \bx_1^{(p)^T} & a^{(p)}
\endpmatrix,
\end{equation}
where $b$ is the smallest rectangular block containing all the
corner elements.

A factorization similar to that in \eqref{parfact}--\eqref{parG}
(the obvious differences are related to the first and last (block)
rows) produces a corresponding {\em reduced system} with the {\em
reduced matrix}
\begin{equation} \label{redmatr1}
T_p = \pmatrix{ccccc}
\alpha^{(0)} & \gamma^{(1)} && & \beta^{(0)} \\
\beta^{(1)} & \alpha^{(1)} & \gamma^{(2)} \\
&\beta^{(2)} & \alpha^{(2)} & \ddots \\
&& \ddots & \ddots & \gamma^{(p)} \\
&&& \beta^{(p)} & \alpha^{(p)}
\endpmatrix.
\end{equation}

\begin{figure}[t]
\begin{center}
\includegraphics*[width=9cm]{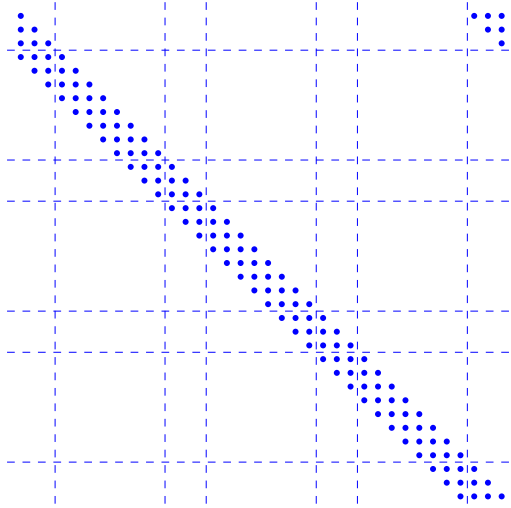}
\vspace*{-1cm}
\end{center}
\caption{\it Partitioning of the matrix in Fig.\,\ref{Fig3b} after
row permutation. Each point represents a (block) entry of the
matrix.} \label{Fig3}
\end{figure}

We observe that, for the very important classes of BABD and
circulant-like matrices (the latter, after a suitable row
permutation, see Figure \ref{Fig3}), both the matrix
\eqref{parpart1} and the {\em reduced matrix} \eqref{redmatr1}
have the form of a lower block bidiagonal matrix (i.e.,
$\cx_j^{(i)}=0$ and $\gamma^{(i)}=0$ for all $i$ and $j$) with an
additional right-upper corner block:
$$A = \pmatrix{ccccccc}
a^{(0)} &  &&&&& b \\
\bx_0^{(1)} & A^{(1)}  \\
&\bx_1^{(1)^T} & a^{(1)}  \\
&& \bx_0^{(2)} & \ddots \\
&&& \ddots & a^{(p-1)}  \\
&&&& \bx_0^{(p)} & A^{(p)}  \\
&&&&& \bx_1^{(p)^T} & a^{(p)}
\endpmatrix,
$$ and
$$
T_p = \pmatrix{cccc}
\alpha^{(0)} & & & \beta^{(0)} \\
\beta^{(1)} & \alpha^{(1)}  \\
& \ddots & \ddots  \\
&& \beta^{(p)} & \alpha^{(p)}
\endpmatrix.
$$

We have also to note that, for this kind of matrices, the overall
computational cost of a parallel factorization algorithm has a
very small increase. On a sequential machine, supposing to
maintain the same partitioning of the matrix on $p>1$ processors,
we have a computational cost which is similar to that of any
efficient sequential algorithm (the corner block $b$ implies the
construction of a fill-in (block) vector) but with better
stability properties (see \cite{Wr}). For this reason, a method
that is widely and efficiently applied to matrices in the form
\eqref{parpart1}, also on a sequential computer, is cyclic
reduction (see \cite{AmGlRo1,AmPa3,AmRo} where cyclic reduction is
applied to BABD matrices).

\section{Parallel solution of differential equations}\label{par}

The numerical solution of ODE-BVPs leads to the solution of large
and sparse linear systems of equations that represent the most
expensive part of a BVP code. The sparsity structure of the
obtained problem depends on the methods implemented. In general,
one-step methods lead to ABD or BABD matrices, depending on the
boundary conditions (separated or not, respectively), while
multistep methods lead to block banded systems (with additional
corner blocks in case of non-separated boundary conditions)
\cite{AmBr97a,AmBr97,BrTr98}. The parallel algorithms previously
described perfectly cope with this kind of systems. For this
reason we do not investigate further on ODE-BVPs.

Conversely, we shall now consider the application of {\em parallel
factorizations} for deriving parallel algorithms for numerically
solving ODE-IVPs, which we assume, for sake of simplicity, to be
linear and in the form

\begin{equation}\label{prob}
y' = Ly+g(t), \qquad t\in[t_0,T], \qquad y(t_0)=y_0\in\RR^m,
\end{equation}
which is, however, sufficient to grasp the main features of the
approach
\cite{AmBr97a,AmBr97,AmBr98,AmBr09,BrTr97,BrTr98a,BrTr98}.

Let us consider a suitable {\em coarse mesh}, defined by the
following partition of the integration interval in \eqref{prob}:

\begin{equation}\label{coarse}
t_0 \equiv \tau_0 < \tau_1 < \cdots < \tau_p \equiv
T.\end{equation}

Suppose, for simplicity, that inside each sub-interval we apply a
given method with a constant stepsize

\begin{equation}\label{hi} h_i = \frac{\tau_i-\tau_{i-1}}N, \qquad
i = 1,\dots,p,\end{equation} to approximate the problem

\begin{equation}\label{probi}
y' = Ly +g(t), \quad t\in[\tau_{i-1},\tau_i], \qquad y(\tau_{i-1})
= y_{0i}, \qquad i = 1,\dots,p.\mbox{~}\end{equation} If $y(t)$
denotes the solution of problem \eqref{prob}, and we denote by

\begin{equation}\label{yni}
y_{ni} \approx y(\tau_{i-1}+nh_i), \qquad n=0,\dots,N, \quad
i=1,\dots,p,\end{equation} the entries of the discrete
approximation, then, in order for the numerical solutions of
\eqref{prob} and \eqref{probi} to be equivalent, we require that
(see \eqref{coarse} and \eqref{yni})

\begin{equation}\label{initi}
y_{01} = y_0, \qquad y_{0i}\equiv y_{N,i-1}, \quad
i=2,\dots,p.\end{equation} For convention, we also set
\begin{equation}\label{init0} y_{01}\equiv y_{N0}.\end{equation}

Let now suppose that the numerical approximations to the solutions
of \eqref{probi} are obtained by solving discrete problems in the
form

\begin{equation}\label{M_iy_i}
M_i \yx_i = \vx_i y_{0i} +\gx_i, \qquad \yx_i =
\pmatrix{c}y_{1i},~\dots~,~y_{Ni}\endpmatrix^T,\qquad
i=1,\dots,p,\end{equation} where the matrices $M_i\in\RR^{mN\times
mN}$ and $\vx_i\in\RR^{mN\times m}$, and the vector
$\gx_i\in\RR^{mN}$, do depend on the chosen method (see, e.g.,
\cite{AmBr97a,AmBr97}, for the case of block BVMs) and on the
problems \eqref{probi}. Clearly, this is a quite general
framework, which encompasses most of the currently available
methods for solving ODE-IVPs. By taking into account all the above
facts, one obtains that the global approximation to the solution
of \eqref{prob} is obtained by solving a discrete problem in the
form (hereafter, $I_r$ will denote the identity matrix of
dimension $r$):

\begin{eqnarray}\nonumber
M\yx \equiv\pmatrix{rrrrr}
I_m\\
-\vx_1 & M_1\\
       & -V_2 &M_2\\
       & &\ddots&\ddots\\
       & &      &-V_p &M_p\endpmatrix
     \pmatrix{c} y_{N0} \\\yx_1\\ \yx_2\\ \vdots\\ \yx_p\endpmatrix
       &=&\pmatrix{c} y_0\\\gx_1\\ \gx_2\\ \vdots\\ \gx_p\endpmatrix,
       \\ \label{My}\\ \nonumber
       V_i = [O\,|\, \vx_i]\in\RR^{mN\times mN}, &&\qquad i=2,\dots,p.
       \end{eqnarray}

Obviously, this problem may be solved in a sequential fashion, by
means of the iteration (see \eqref{initi}-\eqref{init0}):

$$y_{N0} = y_0, \qquad
M_i\yx_i =\gx_i+\vx_i y_{N,i-1},\qquad i=1,\dots,p.$$

Nevertheless, by following arguments similar to those in the
previous section, we consider the factorization:

$$M =
\pmatrix{ccccc} I_m \\&M_1\\ &&M_2\\ &&&\ddots \\&&&&M_p
\endpmatrix
\pmatrix{rrrrrr}
I_m\\
-\wx_1&~I_{mN}\\
      &-W_2 &~I_{mN}\\
      &    &\ddots&\ddots\\
      &    &      &-W_p &~I_{mN}\endpmatrix,$$
where (see \eqref{My})

\begin{equation}\label{Wi}
W_i = [O\,|\, \wx_i]\in\RR^{mN\times mN}, \qquad \wx_i =
M_i^{-1}\vx_i\in\RR^{mN\times m}.
\end{equation}
Consequently, at first we solve, in parallel, the systems

\begin{equation}\label{D}
M_i\zx_i= \gx_i, \qquad \zx_i =
\pmatrix{c}z_{1i},~\dots~,~z_{Ni}\endpmatrix^T,\qquad i=1,\dots,p,
\end{equation}
and, then, (see \eqref{Wi} and \eqref{initi}) recursively update the
local solutions,

\begin{eqnarray}\nonumber
\yx_1 &=& \zx_1 +\wx_1 y_{01},\\
\label{W}\\ \nonumber \yx_i &=& \zx_i+W_i\yx_{i-1} \equiv \zx_i
+\wx_i y_{0i}, \qquad i=2,\dots,p.
\end{eqnarray}

The latter recursion, however, has still much parallelism. Indeed,
if we consider the partitionings (see \eqref{M_iy_i}, \eqref{D}, and
\eqref{Wi})

\begin{equation}\label{part}
\yx_i= \pmatrix{c} \hbfy_i\\ y_{Ni}\endpmatrix, \qquad \zx_i =
\pmatrix{c} \hbfz_i\\ z_{Ni}\endpmatrix, \qquad \wx_i =
\pmatrix{c}\hbfw_i\\ w_{Ni}\endpmatrix, \quad w_{Ni}\in\RR^{m\times
m},\end{equation} then \eqref{W} is equivalent to solve, at first,
the {\em reduced system}

\begin{equation}\label{redsys}
\pmatrix{cccc}
I_m\\
-w_{N1}&I_m\\
      &\ddots&\ddots\\
      &      &-w_{N,p-1}&I_m\endpmatrix
      \pmatrix{c} y_{01}\\y_{02}\\ \vdots\\ y_{0p}\endpmatrix =
      \pmatrix{c} y_0\\ z_{N1}\\ \vdots \\ z_{N,p-1}\endpmatrix,
      \end{equation}
i.e.,
\begin{equation}\label{recur}
y_{01} = y_0, \qquad y_{0,i+1} = z_{Ni}+w_{Ni}y_{0i}, \qquad
i=1,\dots,p-1,
\end{equation}
after which performing the $p$ parallel updates

\begin{equation}\label{pups}
\hbfy_i = \hbfz_i +\hbfw_i y_{0i}, \qquad i=1,\dots,p-1, \qquad
\yx_p = \zx_p+\wx_p y_{0p}.
\end{equation}
We observe that:
\begin{itemize}

\medskip
\item the parallel solution of the  $p$ systems in \eqref{D} is
equivalent to compute the approximate solution of the following
$p$ ODE-IVPs,
\begin{equation}\label{probi0}
z' = Lz +g(t), \quad t\in[\tau_{i-1},\tau_i], \quad z(\tau_{i-1}) =
0, \quad i = 1,\dots,p,\end{equation} in place of the corresponding
ones in \eqref{probi};

\medskip
\item the solution of the {\em reduced system} \eqref{redsys}-\eqref{recur} consists
in computing the proper initial values $\{y_{0i}\}$ for the
previous ODE-IVPs;

\medskip
\item the parallel updates \eqref{pups} update the approximate
solutions of the ODE-IVPs \eqref{probi0} to those of the
corresponding ODE-IVPs in \eqref{probi}.

\end{itemize}

\medskip
\begin{rem}\label{uno} Clearly, the solution of the first
(parallel) system in \eqref{D} and the first (parallel) update in
(\ref{W}) (see also \eqref{pups}) can be executed together, by
solving the linear system (see \eqref{initi})
\begin{equation}\label{y1} M_1\yx_1 =
\gx_1+\vx_1y_0,\end{equation} thus directly providing the final
discrete approximation on the first processor; indeed, this is
possible, since the initial condition $y_0$ is given.
\end{rem}\medskip

We end this section by emphasizing that one obtains an almost
perfect parallel speed-up, if $p$ processors are used, provided
that the cost for the solution of the {\em reduced system}
\eqref{redsys} and of the parallel updates \eqref{pups} is small,
with respect to that of \eqref{D}  (see \cite{AmBr97a,AmBr97} for
more details). This is, indeed, the case when the parameter $N$ in
\eqref{hi} is large enough and the coarse partition \eqref{coarse}
can be supposed to be {\em a priori} given.

\section{Connections with the ``Parareal'' algorithm}\label{parareal}

We now briefly describe the ``Parareal'' algorithm introduced in
\cite{LiMaTu01,MaTu02}, showing the existing connections with the
parallel method previously described. This method, originally
defined for solving PDE problems, for example linear or
quasi-linear parabolic problems, can be directly cast into the ODE
setting via the semi-discretization of the space variables; that
is, by using the method of lines. In more detail, let us consider
the problem

\begin{equation}\label{pde}
\ddt y = \cL\, y, \qquad t\in[t_0,T], \qquad y(t_0)=y_0,
\end{equation}
where $\cL$ is an operator from a Hilbert space $V$ into $V'$. Let
us consider again the partition \eqref{coarse} of the time interval,
and consider the problems

\begin{equation}\label{pdei}
\ddt y = \cL\, y, \qquad t\in[\tau_{i-1},\tau_i], \qquad
y(\tau_{i-1})=y_{0i}, \qquad i=1,\dots,p.
\end{equation}
Clearly, in order for \eqref{pde} and \eqref{pdei} to be
equivalent, one must require that

\begin{equation}\label{initipde}
y_{0i}=y(\tau_{i-1}), \qquad i=1,\dots,p.
\end{equation}
The initial data \eqref{initipde} are then formally related by means
of suitable {\em propagators} $\cF_i$ such that

\begin{equation}\label{initipde1}
y_{0,i+1} = \cF_i y_{0i}, \qquad i=1,\dots,p-1.
\end{equation}

The previous relations can be cast in matrix form as ($\cI$ now
denotes the identity operator)

\begin{equation}\label{redsyspde1}
F \yx \equiv \pmatrix{cccc}
~\cI\\
-\cF_1 &~\cI\\
       &\ddots &\ddots\\
       &       &-\cF_{p-1} &~\cI\endpmatrix
       \pmatrix{c} y_{01} \\ y_{02}\\ \vdots \\ y_{0p}\endpmatrix
       = \pmatrix{c} y_0\\ 0\\ \vdots\\ 0\endpmatrix \equiv \bfeta.
       \end{equation}
For solving \eqref{redsyspde1}, the authors essentially define the
splitting

$$F = (F-G)+G, \qquad G = \pmatrix{cccc}
~\cI\\
-\cG_1 &~\cI\\
       &\ddots &\ddots\\
       &       &-\cG_{p-1} &~\cI\endpmatrix,
$$
with {\em coarse propagators}

$$\cG_i\approx \cF_i, \qquad i=1,\dots,p,$$
and consider the iterative procedure

$$G \yx^{(k+1)} = (G-F)\yx^{(k)} +\bfeta, \qquad k=0,1,\dots,$$
with an obvious meaning of the upper index. This is equivalent to
solve the problems

\begin{eqnarray}\nonumber
y_{01}^{(k+1)} &=& y_0,\\
\label{parareal1} y_{0,i+1}^{(k+1)} &=& \cG_i y_{0i}^{(k+1)} +
(\cF_i-\cG_i)y_{0i}^{(k)}, \qquad i=1,\dots,p-1,
\end{eqnarray}
thus providing good parallel features, if we can assume that the
coarse operators $\cG_i$ are ``cheap'' enough. The iteration
\eqref{parareal1} defines the ``Parareal'' algorithm, which is
iterated until

$$\|y_{0i}^{(k+1)}-y_{0i}^{(k)}\|, \qquad i=2,\dots,p,$$
are suitably small. In the practice, in case of linear operators,
problem \eqref{pde} becomes, via the method of lines, an ODE in
the form \eqref{prob}, with $L$ a huge and very sparse matrix.
Consequently, problems \eqref{pdei} become in the form
\eqref{probi}. Similarly, the propagator $\cF_i$ consists in the
application of a suitable discrete method for approximating the
solution of the corresponding $i$th problem in \eqref{probi}, and
the coarse propagator $\cG_i$ describes the application of a much
cheaper method for solving the same problem. As a consequence, if
the discrete problems corresponding to the propagators $\{\cF_i\}$
are in the form \eqref{M_iy_i}, then the discrete version of the
recurrence \eqref{initipde1} becomes exactly \eqref{recur}, as
well as the discrete counterpart of the matrix form
\eqref{redsyspde1} becomes \eqref{redsys}.

\medskip We can then conclude that the ``Parareal'' algorithm in
\cite{LiMaTu01,MaTu02} {\em exactly} coincides with the iterative
solution of the {\em reduced system} \eqref{redsys}, induced by a
suitable splitting.
\smallskip

We observe that the previous iterative procedure may be very
appropriate, when the matrix $L$ is large and sparse since, in
this case, the computations of the block vectors $\{\wx_i\}$ in
\eqref{Wi}, and then of the matrices $\{w_{Ni}\}$ (see
\eqref{part}) would be clearly impractical. Moreover, it can be
considerably improved by observing that
$$w_{Ni}y_{0i} \approx {\rm e}^{(\tau_i-\tau_{i-1})L}y_{0i}.$$
Consequently, by considering a suitable approximation to the
matrix exponential, the corresponding parallel algorithm turns out
to become semi-iterative and potentially very effective, as
recently shown in \cite{AmBr09}.

\end{document}